\documentclass
{article}
\usepackage[T2A]{fontenc}
\usepackage[cp1251]{inputenc}
\usepackage[english,russian]{babel}
\usepackage{amsmath}
\usepackage{amsfonts,amssymb,mathrsfs,amscd}
\usepackage{verbatim}
\usepackage{eucal}
\usepackage{graphicx}
\usepackage{enumerate}
\usepackage{esint}
\usepackage{upgreek}
\let\No=\textnumero

\usepackage[hyper]{izv2e}
\JournalName{}

\numberwithin{equation}{section}

\theoremstyle{plain}
\newtheorem{theorem}{Теорема}

\newtheorem{theo}{Теорема} 

\newtheorem{lemma}{Лемма}
\newtheorem{propos}{Предложение}

\theoremstyle{definition}
\newtheorem{proof}{Доказательство}
\newtheorem{remark}{Замечание}

\renewcommand{\leq}{\leqslant} 
\renewcommand{\geq}{\geqslant}
\newcommand{\RR}{\mathbb{R}} 
\newcommand{\CC}{\mathbb{C}} 
\newcommand{\NN}{\mathbb{N}} 

\newcommand{\DD}{\mathbb{D}} 

\DeclareMathOperator{\rad}{rad}
\DeclareMathOperator{\bal}{bal}
\DeclareMathOperator{\Bal}{Bal}

\DeclareMathOperator{\supp}{supp} 
\DeclareMathOperator{\type}{type} 
\DeclareMathOperator{\strip}{str} 
\DeclareMathOperator{\ord}{ord}
\DeclareMathOperator{\rh}{rh}
 \DeclareMathOperator{\lh}{lh}
 \DeclareMathOperator{\up}{up}
\DeclareMathOperator{\dd}{d}
\renewcommand{\Re}{\operatorname{Re}}
\renewcommand{\Im}{\operatorname{Im}}

\begin{document} 

\title{Выметание распределений зарядов и субгармонических функций на полосу
}
\author[B.\,N.~Khabibullin]{Б.\,Н.~Хабибуллин}
\address{Башкирский государственный университет}
\email{khabib-bulat@mail.ru} 

\date{12.04.2022}
\udk{517.574}

\maketitle
	
	\begin{fulltext}

\begin{abstract} 
Рассматриваются две конструкции выметания на комплексной плоскости $\mathbb C$ с вещественной осью $\mathbb R$ при  $0\leq b\in \mathbb R$. 

Пусть $u\not\equiv -\infty$ --- субгармоническая  функций на $\mathbb C$ порядка 
$$
\operatorname{ord}[u]:=\limsup_{z\to \infty} \frac{\ln \max\{1,u(z)\}}{\ln |z|}\leq 1,
$$  
$U=u-v$ --- разность субгармонических  функций $u$ и  $v\not\equiv -\infty$ с   $\operatorname{ord}[v]\leq 1$, на $\mathbb C$, т.е. $\delta$-субгармоническая функция на $\mathbb C$ порядка $\operatorname{ord}[U]\leq 1$. 
Тогда  найдётся  $\delta$-субгармоническая функция $V\not\equiv \pm\infty$ на 
$\mathbb C$ порядка $\operatorname{ord}[V]\leq 1$, гармоническая  на $\bigl\{ z \in \mathbb C\bigm| |\Re z|> b\bigr\}$, а также полярное множество $E\subset \mathbb C$, для которых $U(z)\equiv V(z)$ при всех $z\in \bigl\{ z \in \mathbb C\bigm| |\Re z|\leq b\bigr\}\setminus E$.

Если  $u$ --- субгармоническая функция конечного типа при порядке $1$, т.е. $\limsup_{z\to \infty} u(z)/|z|<+\infty$, то существуют субгармонические функции   $u_{\mathbb R}$ и $u_b$  конечного типа при порядке $1$, гармонические соответственно 
на $\mathbb C\setminus \mathbb R$ и $\bigl\{ z \in \mathbb C\bigm| |\Re z|> b\bigr\}$ , для которых одновременно 
$$
\begin{cases}
u(z)\equiv u_{\mathbb R}(z)+u_b(z) \text{ при всех $z\in {\mathbb R}\bigcup \bigl\{ z \in \mathbb C\bigm| |\Re z|\leq b\bigr\}$},\\ 
u(z)\leq u_{\mathbb R}(z) +u_b(z)  \text{ при всех $z\in \mathbb C$.}
\end{cases}
$$
При этом прослеживаются взаимосвязи между различными логарифмическими характеристиками распределений масс и зарядов Рисса 
субгармонических и $\delta$-суб\-г\-а\-р\-м\-о\-н\-и\-ч\-е\-с\-ких функций.

Библиография:  14  наименований

\end{abstract}
	

\markright{Выметание распределений зарядов и субгармонических функций}

\footnotetext[0]{Исследование выполнено за счёт гранта Российского научного фонда №~22-21-00026, 
\href{https://rscf.ru/project/22-21-00026/}{https://rscf.ru/project/22-21-00026/}}

Основные результаты статьи --- теоремы \ref{Balvs} и \ref{prop4_1} соответственно из \S\S~\ref{S3},\ref{bRu}. 

Можно сразу перейти  к \S~\ref{BalCrh}, а к \S~\ref{11def}  обращаться по мере необходимости.

\section{Некоторые  обозначения, определения, соглашения}\label{11def} 
 Одноточечные множества $\{a\}$ часто записываем без фигурных скобок, т.е. просто как $a$. Так,   $\NN_0:=0\cup \NN=\{0,1, \dots\}$ для  множества $\mathbb N:=\{1,2, \dots\}$ {\it натуральных чисел,\/}  
$\CC_{\infty}:=\CC\cup \infty$ и $\overline \RR:=-\infty \cup \RR\cup +\infty$ --- {\it расширенные\/}
комплексная плоскость и вещественная ось с $-\infty:=\inf \RR\notin \RR$, $+\infty:=\sup \RR\notin \RR$, неравенствами $-\infty\leq x\leq +\infty$ для любого $x\in \overline \RR$
 и естественной порядковой топологией. По определению 
$\sup \varnothing:=-\infty$ и $\inf \varnothing:=+\infty$ для {\it пустого множества\/} $\varnothing$.
Символом  $0$, кроме нуля,   могут обозначаться {\it нулевые\/} функции,  меры  и пр.
Для $x\in X\subset \overline \RR$ его {\it положительную часть\/} обозначаем как $x^+:=\sup\{0,x \}$, $X^+:=\bigl\{x^+\bigm|  x\in X\bigr\}$. {\it Расширенной числовой функции\/} $f\colon S\to \overline \RR$ сопоставляем её {\it положительную часть\/} $f^+\colon s\underset{s \in S}{\longmapsto} (f(s))^+\in \overline{\RR}^+$ 
и {\it отрицательную часть\/} $f^-:=(-f)^+\colon S\to \overline{\RR}^+ $. Как обычно, пишем $f\not\equiv c$, если функция $f$ принимает хотя бы одно значение, отличное от $c$, в области её определения. 

Для $x_0\in \RR$ и расширенной числовой функции $m\colon x_0+\RR^+ \to \overline \RR$ определим  
\begin{equation*}
\ord[m]:=\limsup_{x\to +\infty} \frac{\ln \bigl(1+m^+(x)\bigr)}{\ln x}\in
\overline \RR^+
 \end{equation*} 
{\it --- порядок\/}  (роста) функции $m$ (около $+\infty$), а  для $p\in \RR^+$
\begin{equation*} 
\type_p[m]:=\limsup_{x\to +\infty} \frac{m^+(x)}{x^p}\in \overline \RR^+
\label{typevf}
 \end{equation*} 
{\it --- тип\/} (роста) функции $m$ {\it при порядке\/} $p$ (около $+\infty$) \cite{Boas}, \cite{Levin56}, \cite{Levin96}, 
\cite{Kiselman}, \cite[2.1]{KhaShm19}, а для произвольной 
 функции $u\colon \CC\to \overline \RR$ с  {\it радиальный функцией роста\/}
\begin{equation*} 
{\mathrm M}_u \colon r\underset{r\in \RR^+}{\longmapsto}
\sup\bigl\{u(z)\bigm| |z|=r\bigr\}
\label{u}
\end{equation*} 
по определению $\ord[u]:=\ord[{\mathrm M}_u]$ и $\type_p[u]:=\type_p[{\mathrm M}_u]$ 
--- соответственно {\it порядок\/} и 
{\it тип функции $u$ при порядке $p$} \cite{Boas}, \cite{Levin96}, \cite{Kiselman}, \cite[Замечание 2.1]{KhaShm19}.  Функции $u$ конечного типа 
$\type_1[u]\in \RR^+$ при порядке $p=1$ называем просто функциями {\it конечного типа,\/} не упоминая  и не указывая 
порядок $1$ в $\type[u]:=\type_1[u]$. 

{\it Распределением масс\/} --- это  {\it положительная   мера Радона\/}  \cite{EG},  \cite[Дополнение A]{Rans}, \cite[гл.~3]{HK}, 
а {\it распределение зарядов\/} ---  разность распределений масс \cite{Landkof}.
Для распределений масс   или зарядов  {\it на\/} $\CC$, как правило, не указываем, где  они заданы.   
Для  субгармонической  в области из $\CC$ функции $u\not\equiv -\infty$ действие  на неё   {\it оператора Лапласа\/} ${\bigtriangleup}$  в смысле теории обобщённых функций  определяет её {\it распределение масс Рисса\/}
\begin{equation}\label{Riesz}
\frac{1}{2\pi}{\bigtriangleup}u=:\varDelta_u
\end{equation}
 в этой области  \cite{HK}, \cite{Rans}, \cite{Az}. Для обозначения распределения масс Рисса функции $u$ используем  как первую  форму записи $\frac{1}{2\pi}{\bigtriangleup}u$ из \eqref{Riesz}, так и вторую  $\varDelta_u$.
 Далее $D_z(r):=\bigl\{w \in \CC \bigm| |w-z|<r\bigr\}$ и  $\overline{D}_z(r):=\bigl\{w \in \CC_{\infty} \bigm| |w-z|\leq r\bigr\}$,
а также $\partial \overline{D}_z(r):=\overline{D}_z(r)\setminus {D}_z(r)$
---  соответственно {\it открытый} и   {\it замкнутый круги,\/} а также  {\it окружность  радиуса\/ $r\in \overline \RR^+$
 с центром\/} $z\in \CC$, а   $\DD:=D_0(1)$ и  $\overline \DD:=\overline D_0(1)$, а также  
$\partial \overline \DD:=\partial \overline{D}_0(1)$ --- соответственно {\it открытый\/} и {\it замкнутый 
единичные круги,\/} а также {\it единичная окружность} в  $\CC$. 

Через $\CC_{ \rh}:=\bigl\{z\in \CC \bigm| \Re z>0\bigr\}$ и $\CC_{\overline  \rh}:=\CC_{ \rh}\cup i\RR$, а также
$\CC_{ \lh}:=-\CC_{ \rh}$ и $\CC_{\overline \lh}:=-\CC_{\overline  \rh}$ обозначаем 
соответственно {\it правые открытую} и   {\it замкнутую полуплоскости,\/} а также  
{\it левые открытую} и   {\it замкнутую  полуплоскости\/} в  $\CC$. 

Для распределения зарядов  $\nu$ на $S\subset \CC$ через $\nu^+:=\sup\{\nu,0\}$, $\nu^-:=(-\nu)^+$
и $|\nu|:=\nu^++\nu^-$ обозначаем соответственно {\it верхнюю, нижнюю\/}
и {\it полную вариации\/} распределения зарядов  $\nu$, а $\supp \nu=\supp |\nu|$ --- его {\it носитель,} но  распределение зарядов  $\nu$ {\it сосредоточено на $\nu$-измеримом подмножестве $S_0\subset S$,\/} если полная вариация $|\nu|$ дополнения $S\setminus S_0$ множества $S$ равна нулю. 

{\it Сужение\/} функции $f$ на  $S\subset \CC$ обозначаем как   $f{\lfloor}_S$. Аналогично через $\nu{\lfloor}_S$
обозначается  обозначается и   {\it сужение\/}  положительной меры Бореля или  распределения зарядов $\nu$ на  $\nu$-измеримое  $S\subset \CC$.
При $r\in \overline \RR^+$ для таких   $\nu$  через 
\begin{equation}\label{df:nup} 
\nu_z^{\rad} (r):=\nu \bigl(\,\overline D_z(r)\bigr),\quad \nu^{\rad}(r):=\nu_0^{\rad}(r)
=\nu\bigl(r\overline \DD\bigr)
\end{equation} 
обозначаем  {\it радиальные непрерывные справа считающие  функции  распределения зарядов  $\nu$ с
центрами\/}  соответственно {\it в точке\/ $z\in \CC$} и {\it в нуле.\/}

{\it Верхняя плотность распределения зарядов\/ $\nu$ при порядке\/} $p\in \RR^+$  равна 
\begin{equation}\label{typenu}
\type_p[\nu]:=\type_p\bigl[|\nu|\bigr]
\overset{\eqref{typevf}}{:=} \limsup_{0<r\to +\infty} \frac{|\nu|(r\overline \DD)}{r^p}
\overset{\eqref{df:nup}}{=} \limsup_{0<r\to +\infty} \frac{|\nu|^{\rad}(r)}{r^p}
\in \overline \RR^+,
\end{equation} 
и при $p=1$ упоминание о порядке  опускаем.   В частности, {\it  распределение зарядов\/} $\nu$
{\it конечной верхней плотности,\/} если $\type[\nu]:=\type_1[\nu]<+\infty$. {\it Порядок распределения зарядов\/}
$\nu$ определяется как $\ord[\nu]:=\ord\bigl[|\nu|^{\rad}\bigr]$
через \eqref{df:nup}. 

Для  распределения зарядов $\nu$   
\begin{align}
\ell_{\nu}^{\rh}(r, R)&
:=\int_{r < | z|\leq R} \Re^+ \frac{1}{ z} \dd \nu(z)\in \RR, \quad 0< r < R < +\infty ,
\label{df:dDlm+}\\
\ell_{\nu}^{\lh}(r, R)&
:=\int_{r< |z|\leq R}\Re^- \frac{1}{z} \dd \nu(z)\in \RR,  \quad 0< r < R <  +\infty ,
\label{df:dDlm-}
\end{align}
---  соответственно {\it правая\/} и {\it левая логарифмические функции интервалов\/} $(r,R]$ на $\RR^+$.  В случае  {\it распределения масс\/} $\mu$ они порождают его 
{\it правую\/} $\ell_{\mu}^{\rh}$ и {\it левую $\ell_{\mu}^{\lh}$ логарифмические меры\/} на $\overline \RR^+\setminus 0$, допуская  и $R=+\infty$  в \eqref{df:dDlm+}--\eqref{df:dDlm-} с возможными  значениями $+\infty$ для $\ell_{\mu}^{\rh}(r, +\infty)$ 
и $\ell_{\mu}^{\lh}(r, +\infty)$,  а также его двустороннюю  {\it логарифмическую субмеру} на $\RR^+\setminus 0$, определённую как  
\begin{equation}\label{df:dDlLm}
\ell_{\mu}(r, R):=\max \bigl\{ \ell_{\mu}^{\lh}(r, R), \ell_{\mu}^{\rh}(r,R)\bigr\}\in \overline \RR^+, \quad 0< r < R \leq +\infty. 
\end{equation}

Для  $b\in \RR^+$  через 
\begin{equation}\label{{strip}c}
 \strip_b:=\Bigl\{z\in \CC\Bigm| |\Re z|< b\Bigr\}, \quad 
\overline \strip_b:=\Bigl\{z\in \CC\Bigm| |\Re z|\leq b\Bigr\}
\end{equation}
обозначаем {\it вертикальные\/} соответственно {\it открытую\/} и {\it  замкнутую полосу 
ширины\/ $2b$ со средней линией $i\RR$.}

\section{Выметание  из  правой полуплоскости}\label{BalCrh}

\subsection{Выметания рода $0$ и $1$ распределений зарядов}
Для распределения заряда $\nu$ используем \cite[формула (1.9)]{KhaShm19}   его {\it функцию  распределения на\/}  $\RR$, 
обозначаемую как  $\nu_{\RR}\colon \RR\to  \RR$ и определённую  равенствами 
\begin{equation}\label{nuR} 
\nu_{\RR}(x_2)-\nu_{\RR}(x_1):=\nu\bigl((x_1,x_2]\bigr), \quad
-\infty <x_1<x_2<+\infty, 
\end{equation}
и   {\it функцию распределения $\nu_{i\RR}\colon \RR\to \RR$    на\/} $i\RR$, определённую равенствами
\begin{equation}\label{nuiR} 
\nu_{i\RR}(y_2)-\nu_{i\RR}(y_1):=\nu\bigl(i(y_1,y_2]\bigr), \quad -\infty <y_1<y_2<+\infty. 
\end{equation}
Поскольку эти функции распределения определены лишь с точностью до аддитивной постоянной, 
при необходимости используем их {\it нормировки в нуле}
\begin{equation}\label{nuo}
\nu_{\RR}(0):=0, \quad \nu_{i\RR}(0):=0.
\end{equation}
По построению \eqref{nuR} и \eqref{nuiR}    функции  $\nu_{\RR}$ и $\nu_{i\RR}$  локально ограниченной вариации  на $\RR$, а в случае {\it распределения масс\/} $\nu$ обе эти функции возрастающие. 
Обратно, любая функция локально ограниченной вариации на $\RR$ или $i\RR$ однозначно определяет распределение зарядов с носителем соответственно на  $\RR$ или $i\RR$. 

Мы напоминаем и адаптируем основные понятия и утверждения из \cite{KhaShm19} и \cite{KhaShmAbd20},
 а также частично из \cite{Kha91} и \cite{Kha91AA} о  выметании конечного рода $q\in \NN_0$ распределений зарядов, но пока применительно только к правой полуплоскости $\CC_{\rh}$ в случаях  $q:=0$ и $q:=1$. В 
\cite{KhaShm19} и \cite{KhaShmAbd20} в основном рассматривается {\it верхняя полуплоскость\/} $\CC^{\up}:=i\CC_{ \rh}$, что переносится на $\CC_{\rh}$ поворотом на прямой угол.

{\it Характеристическую функцию множества\/}  $S$ обозначаем через 
\begin{equation}\label{SdrS}
\boldsymbol{1}_S\colon z\underset{z\in \mathbb C}{\longmapsto} \begin{cases}
1&\text{ если $z\in S$},\\
0&\text{ если $z\notin S$}.
\end{cases}
\end{equation}
 
{\it Гармоническая  мера  для\/ $\CC_{\rh}$ в точке\/ $z\in \CC_{\rh}$} на интервалах $i(y_1,y_2]\subset i\overline\RR$
\begin{equation}\label{omega}
\omega_{\rh} \bigl(z,i(y_1,y_2]\bigr){\overset{\text{\cite[3.1]{KhaShm19}}}{:=}}\omega_{\CC_{\rh}}(z,i(y_1,y_2])
\underset{z\in \CC_{\rh}}{:=}\frac1{\pi}
\int_{y_1}^{y_2}\Re \frac{1}{z-iy} \dd y 
\end{equation} 
равна делённому на $\pi$ углу, под которым виден интервал $i(y_1,y_2]$ из точки $z\in \CC_{\rh}$ \cite[(3.1)]{Kha91AA}, \cite[1.2.1, 3.1]{KhaShm19},  а  в точках  мнимой оси $i\RR$ определяется как 
\begin{equation}\label{oiR}
\omega_{\rh} \bigl(iy,i(y_1,y_2]\bigr):=\boldsymbol{1}_{(y_1,y_2]}(y)
\quad\text{при $y\in  \RR$.}
\end{equation}
Для  распределения зарядов  $\nu $ при {\it классическом условии Бляшке\/} для  $\CC_{\rh}$
\begin{equation}\label{Blcl}
l_{|\nu|}^{\rh}(1,+\infty)\overset{\eqref{df:dDlLm}}{=} \int_{\CC_{\rh}\setminus  \DD} \Re \frac{1}{z}\dd |\nu| (z)<+\infty
\end{equation}
определено \cite[следствие 4.1, теорема 4]{KhaShm19}  
его классическое выметание из $\CC_{\rh}$ на  $ \CC_{\overline \lh}$ с носителем на  $  \CC_{\overline \lh}$, которое  в более широких рамках  \cite[определение 3.1]{KhaShmAbd20} 
представляет собой {\it выметание рода\/ $0$,\/} обозначавшееся в \cite{KhaShmAbd20} как
$\nu^{\bal[0]}_{\CC_{\overline \lh}}$. 
Здесь используется чуть более компактная  запись  
 $\nu^{\bal^0_{\rh}}:=\nu^{\bal[0]}_{\CC_{\overline \lh}}$. По определению {\it распределение зарядов\/}  $\nu^{\bal^0_{\rh}}$ --- это {\it сумма сужения\/}  $\nu{\lfloor}_{\CC_{\lh}}$ на $\CC_{\lh}$  с {\it распределением зарядов на $i\RR$,} определяемым в обозначениях \eqref{nuiR}  {\it функцией распределения\/}
\begin{equation}\label{mubal}
\nu^{\bal^0_{\rh}}_{i\RR}(y_2)-\nu^{\bal^0_{\rh}}_{i\RR}(y_1)\overset{\eqref{omega},\eqref{oiR}}{:=}
\int\limits_{\CC_{\overline \rh}} \omega_{\rh}\bigl(z, i(y_1,y_2]\bigr) \dd \nu(z)
\end{equation}
с нормировкой вида \eqref{nuo} при необходимости. 
Классическое выметание  рода $0$ не увеличивает полную  меру полной вариации  распределения зарядов, 
поскольку гармоническая мера \eqref{omega} вероятностная и 
\begin{equation}\label{omega1}
\bigl|\nu^{\bal^0_{\rh}}\bigr|(S)\underset{S\subset \CC}{\overset{\eqref{mubal}}{\leq}} |\nu|(S).
\end{equation}  

В \cite[определение 2.1]{KhaShmAbd20} вводилось понятие гармонического   заряда  рода\/ $1$ для верхней полуплоскости \/ $\CC^{\up}$ в точке\/ $z\in \CC^{\up}$,  обозначавшегося  в \cite[формула (2.1)]{KhaShmAbd20} через  $\Omega^{[1]}_{\CC^{\up}}$. Здесь используем 
поворот на прямой угол с переходом от $\CC^{\up}$ к $\CC_{\rh}$ и определим 
 {\it гармонический    заряд рода\/ $1$ для правой  полуплоскости \/ $\CC_{\rh}$} 
 как функцию $\Omega_{\rh}$ ограниченных  интервалов $i(y_1,y_2]\subset i\RR$ по правилу 
\begin{equation}\label{Ocr}
\Omega_{\rh}\bigl(z,i(y_1,y_2]\bigr)
\overset{\eqref{omega},\eqref{oiR}}{:=}\omega_{\rh}\bigl(z,i(y_1,y_2]\bigr)-\frac{y_2-y_1}{\pi}\Re\frac{1}{z}
\quad\text{в $z\in \CC_{\overline \rh}\setminus 0$}.
\end{equation}
Для  распределения зарядов  $\nu$ в  \cite[определение 3.1, теорема 1, замечание 3.3]{KhaShmAbd20} определялось 
{\it выметание\/ $\nu_{\CC_{\overline \lh}}^{\bal[1]}$ рода\/ $1$ распределения зарядов\/  $\nu$ из\/  $\CC_{\rh}$ на \/} $ \CC_{\overline \lh}$
  при $0\notin \supp \nu$. Здесь используется  более компактная запись  для такого выметания $\nu^{\bal^1_{\rh}}:=\nu^{\bal[1]}_{ \CC_{\overline \lh}}$.  По определению {распределение зарядов\/} $\nu^{\bal^1_{\rh}}$ --- это 
{\it сумма  сужения\/} $\nu{\lfloor}_{\CC_{\lh}}$ с {\it  распределением зарядов на\/} $i\RR$, определяемым в обозначениях \eqref{nuiR}  {\it функцией распределения} 
\begin{equation}\label{df:nurh}
\nu^{\bal^1_{\rh}}_{i\RR}(y_2)-\nu^{\bal^1_{\rh}}_{i\RR}(y_1)\overset{\eqref{Ocr}}{=}
\int_{\CC_{\rh}} \Omega_{\rh} \bigl(z, i(y_1,y_2]\bigr)\dd \nu (z)
\end{equation}
с нормировкой вида \eqref{nuo} при необходимости.

Ограничение $0\notin \supp \nu$ для выметания рода $1$ легко преодолевается, если скомбинировать выметание рода $0$ части $\nu$  около нуля с выметанием рода $1$ для оставшейся части $\nu$. Для этого определяем {\it комбинированное выметание рода\/ $01$ 
распределения зарядов\/ $\nu$ из\/ $\CC_{\rh}$ на\/ $\CC_{\overline \lh}$} \cite[замечание 3.3, (3.43), (4.1)]{KhaShmAbd20}
\begin{equation}\label{bal01}
\nu^{\bal_{\rh}^{01}}:=\bigl(\nu{\lfloor}_{r_0\DD}\bigr)^{\bal_{\rh}^{0}}+\bigl(\nu{\lfloor}_{\CC\setminus r_0\DD}\bigr)^{\bal_{\rh}^{1}}
\end{equation}
при каком-нибудь фиксированном радиусе  $r_0\in \RR^+\setminus 0$. 
\begin{remark}\label{rem0nu}
Круг $r_0\DD$  в правой части \eqref{bal01}  можно заменить на любое ограниченное борелевское множество, содержащее полукруг $r_0\DD\setminus \CC_{\overline \lh}$.  Можно  обойтись и  без 
выметания рода $0$, положив $\nu^{\bal_{\rh}^{01}}:=\nu^{\bal_{\rh}^{1}}$, если $0\notin \supp \nu$ или, более общ\'о,
\begin{equation*}
\int_{r_0\DD} \Re^+\frac{1}{z} \dd |\nu|(z)<+\infty, 
\end{equation*}
что, очевидно, выполнено, если для некоторого $r_0\in \RR^+\setminus 0$ имеем 
\begin{equation}\label{nur0+}
|\nu|\bigl(r_0\DD\cap \CC_{ \rh}\bigr)=0.
\end{equation}
\end{remark}
Теперь вопрос существования выметания  $\nu^{\bal_{\rh}^{01}}$  упирается лишь в поведение 
распределения зарядов $\nu$ около бесконечности.

Распределение зарядов $\nu$ принадлежит  {\it классу сходимости при порядке\/} роста  $p\in \NN_0$, если
\cite[определение 4.1]{HK}, \cite[\S~2, 2.1, (2.3)]{KhaShm19} 
\begin{equation}\label{sufc}
\int_1^{+\infty}\frac{|\nu|^{\rad}(t)}{t^{p+1}}\dd t<+\infty. 
\end{equation}

Сочетания специальных случаев из  \cite[теорема 1]{KhaShmAbd20},  \cite[теорема 3, п.~4]{KhaShmAbd20} и в некоторой мере   \cite[теорема 3.1]{Kha91AA} легко дают следующие достаточные условия существования выметание распределения зарядов  рода $01$. 

\begin{theorem}\label{theoB1} Пусть $\nu$  ---  распределение зарядов, для которого сужение $\nu{\lfloor}_{\CC_{ \rh}}$ при\-н\-а\-д\-л\-ежит   классу сходимости при порядке $p\overset{\eqref{sufc}}{=}2$.  
Тогда существует выметание $\nu^{\bal^{01}_{\rh}}$ рода $01$ из $\CC_{\rh}$ на $ \CC_{\overline \lh}$. В частности, если $\ord[\nu]<2$, то  $\nu$ из  класса сходимости при порядке $<2$ и   $\ord[\nu^{\bal^{01}_{\rh}}]\leq \ord[\nu]$.  

\end{theorem}

\subsection{Выметания  разностей  субгармонических  функций}\label{6_2}
Пусть $\mathcal U=u-v$ --- разность субгармонических на $\CC$ функций $u$ и $v$, или $\delta$-суб\-г\-а\-р\-м\-о\-н\-и\-ч\-е\-с\-к\-ая функция, для которой при   $u\not\equiv -\infty$ и $v\not\equiv -\infty$ пишем  ${\mathcal U}\not\equiv \pm\infty$.  Значения такой функции ${\mathcal U}\not\equiv \pm\infty$ определены во всех точках, в которых одна из функций  $u$ или $v$ принимает значение из $\RR$, т.е. вне некоторого полярного множества, а её   {\it распределение зарядов Рисса\/}
\begin{equation}\label{Delta}
\varDelta_{\mathcal U}:=\frac{1}{2\pi}{\bigtriangleup}{\mathcal U}\overset{\eqref{Riesz}}{:=}
\frac{1}{2\pi}{\bigtriangleup}u-\frac{1}{2\pi}{\bigtriangleup}v\overset{\eqref{Riesz}}{=}
\varDelta_u-\varDelta_v
\end{equation} 
--- разность распределений масс Рисса  $u$ и $v$. Следуя  \cite[определение 4.1]{KhaShmAbd20}, 
{\it  $\delta$-субгармоническим выметанием $\delta$-субгармонической функции  ${\mathcal U}\not\equiv \pm\infty$  из\/ $\CC_{\rh}$ на $ \CC_{\overline \lh}$}  называем 
каждую  {\it $\delta$-субгар\-м\-о\-н\-и\-ч\-е\-с\-к\-ую  функцию,\/} обозначаемую как ${\mathcal U}^{{\Bal_{\rh}}}$, которая  равна   функции ${\mathcal U}$ на замкнутой левой полуплоскости\/ $\CC_{\overline \lh}$ вне некоторого полярного множества и одновременно гармоническая  на     $\CC_{\rh}$. 

\begin{theorem}[{\rm  (частный случай \cite[теоремы 6 и 7]{KhaShmAbd20})}]\label{Balv}
 Пусть  $\delta$-суб\-г\-а\-р\-м\-о\-н\-и\-ч\-е\-с\-к\-ая  функция ${\mathcal U}\not\equiv \pm \infty$  c распределением зарядов
 Рисса  \eqref{Delta} конечной верхней плотности представима  в виде разности субгармонических функций не более чем первого  порядка.
 Тогда существует $\delta$-суб\-г\-а\-р\-м\-о\-н\-и\-ч\-е\-ское  выметание\/ $\mathcal U^{\Bal_{\rh}}\not\equiv \pm\infty$ из\/ $\CC_{\rh}$ на $ \CC_{\overline \lh}$
c распределением зарядов Рисса 
\begin{equation}\label{112}
\frac{1}{2\pi}{\bigtriangleup}{\mathcal U}^{\Bal_{\rh}}\overset{\eqref{bal01}}{=}\varDelta_{\mathcal U}^{\bal^{01}_{\rh}}, 
\end{equation}  
представимое вне некоторого полярного множества   в виде разности 
\begin{equation}\label{reprvB}
{\mathcal U}^{\Bal_{\rh}}:=u_+-u_-, \quad u_{\pm}\not\equiv -\infty, \quad \ord[u_{\pm}]\leq 1,
\end{equation} 
 субгармонических функций $u_{\pm}\not\equiv -\infty$. Если при этом функция $\mathcal U$ 
гармоническая в открытом полукруге $r_0\DD\cap \CC_{ \rh}$ при некотором $r_0>0$, то правую часть в \eqref{112} можно заменить на выметание $\varDelta_{\mathcal U}^{\bal^{1}_{\rh}}$ рода $1$.  
\end{theorem}

\begin{proof} По \cite[теорема 6]{KhaShmAbd20} для любой $\delta$-субгармонической функции $\mathcal U$ с распределением зарядов Рисса конечного типа  существует выметание ${\mathcal U}^{\Bal_{\rh}}$ с распределением зарядов Рисса  \eqref{112}, представимое в виде 
 \begin{equation}\label{v+u}
{\mathcal U}^{\Bal_{\rh}}=v_+-u_-+H,  \quad v_+\not\equiv -\infty, \quad u_-\not\equiv -\infty,
\end{equation} 
где $v_+$ и $u_-$--- субгармонические функции порядка не выше $1$ , а $H$ --- гармоническая функция на $\CC$.  
При этом если функция $\mathcal U$ представима в виде разности субгармонических функций  порядка не выше $1$, 
то в заключительной части  \cite[теорема 6]{KhaShmAbd20} отмечено, что в качестве  $H$ можно выбрать  гармонический многочлен степени $\deg H\leq 1$. Таким образом, $u_+:=v_++H\not\equiv -\infty$ --- субгармоническая функция порядка не выше $1$ и из 
\eqref{v+u} получаем \eqref{reprvB}. Возможность замены правой части в \eqref{112} на  $\varDelta_{\mathcal U}^{\bal^{1}_{\rh}}$ следует из замечания \ref{rem0nu} в части \eqref{nur0+}. 
\end{proof}

\section{Выметание на вертикальную полосу}\label{S3}

\subsection{Условия Линделёфа}
Распределение зарядов $\nu$ на $\CC$ удовлетворяет  {\it $\RR$-ус\-л\-о\-в\-ию Линделёфа\/}  (рода $1$), если  
\begin{equation}
\sup_{r\geq 1} \biggl| \int_{1<|z|\leq r}\Re\frac{1}{z}\dd \nu(z)\biggr|<+\infty, 
\label{con:LpZR}
\end{equation}
что по определениям \eqref{df:dDlm+}--\eqref{df:dDlm-} эквивалентно соотношению
\begin{equation}
\sup_{r\geq 1} \bigl| \ell^{\rh}_{\nu}(1,r)-\ell^{\lh}_{\nu}(1,r)\bigr|<+\infty, 
\label{con:LpZRl}
\end{equation}
удовлетворяет {\it $i\RR$-условию Линделёфа\/} (рода $1$), если
\begin{equation}
\sup_{r\geq 1}\biggl| \int_{1<|z|\leq r}\Im\frac{1}{z}\dd \nu(z)
\biggr|<+\infty, 
\label{con:LpZiR}
\end{equation}
и удовлетворяет  {\it условию Линделёфа\/} (рода $1$), если
\begin{equation}
\sup_{r\geq 1}\biggl|\int_{1<|z|\leq r}\frac{1}{z}\dd \nu(z)
\biggr|<+\infty. 
\label{con:LpZ}
\end{equation}

Ключевая роль условий Линделёфа отражает следующая классическая  
\begin{theo}[{Вейерштрасса\,--\,Адамара\,--\,Линделёфа\,--\,Брело (\cite[3, Теорема  12]{Arsove53p},  \cite[4.1, 4.2]{HK}, \cite[2.9.3]{Az}, \cite[6.1]{KhaShm19})}]\label{pr:rep}
Если   $u\not\equiv -\infty$ --- субгармоническая функция конечного типа, то её   распределение масс Рисса  $\frac{1}{2\pi}{\bigtriangleup}u$
конечной верхней плотности и удовлетворяет условию Линделёфа \eqref{con:LpZ}. 

Обратно, если  распределение масс    $\nu$ конечной верхней плотности, то существует  субгармоническая функция $u_{\nu}$ с  
$\frac{1}{2\pi}{\bigtriangleup}u_{\nu}=\nu$  порядка $\ord[u_{\nu}]
\leq 1$, которая при выполнении условия Линделёфа \eqref{con:LpZ} для\/  $\nu$   будет уже  функцией  конечного типа. 
При этом  любая  субгармоническая   функция $u$ с  $\frac{1}{2\pi}{\bigtriangleup}u=\nu$
представляется в виде суммы $u=u_{\nu}+H$, где $H$ --- гармоническая функция на\/ $\CC$, которая при условии\/ $\type_2[u]=0$ является  гармоническим многочленом степени\/  $\deg H\leq 1$, а функция $u$ становится функцией порядка\/  $\ord[u]\leq 1$.
\end{theo}

\subsection{Сдвиги и двустороннее  выметание распределения зарядов}\label{5_1Ss}
{\it Зе\-ркальная симметрия $z\underset{z\in \CC}{\longmapsto} -\Bar z$ относительно мнимой оси\/} позволяет все результаты о выметании рода  $1$ из $\CC_{ \rh}$ на $ \CC_{\overline \lh}$ переформулировать для выметания 
из левой полуплоскости $ \CC_{\lh}$ на $\CC_{\overline \rh}$ с заменой, где необходимо, правой логарифмической  функции интервалов \eqref{df:dDlm+} на  левую \eqref{df:dDlm-},  а также с переобозначением  верхнего индекса $^{\bal^1_{\rh}}$  как  $^{\bal^1_{\lh}}$ при выметании из левой полуплоскости $\CC_{ \lh}$.

Для распределения зарядов $\nu$ и точки $w\in \CC$ через 
$\nu_{\vec{w}}$ обозначаем {\it $w$-сдвиг  распределения зарядов\/ $\nu$,} определяемый 
как 
\begin{equation}\label{vecz}
\nu_{\vec{w}}(K):=\nu(K-w) \quad\text{на компактах $K\subset \CC$.}
\end{equation}

\begin{propos}\label{prb2} Пусть $\nu$ --- распределение зарядов конечной верхней пло\-т\-н\-о\-с\-ти. Тогда  
для любых $w\in \CC$ и   $r_0\in \RR^+\setminus 0$  имеем 
\begin{equation}\label{K1}
\sup_{r\geq r_0} \bigl|\ell^{\rh}_{\nu-\nu_{\vec{w}}}(r_0, r)\bigr|
+\sup_{r\geq r_0} \bigl|\ell^{\lh}_{\nu-\nu_{\vec{w}}}(r_0, r)\bigr|<+\infty,
\end{equation}
a  $\nu$ и $\nu_{\vec{w}}$ могут удовлетворять какому-либо одному из трёх видов 
условий Линделёфа \eqref{con:LpZR}--\eqref{con:LpZRl}, \eqref{con:LpZiR} или \eqref{con:LpZ}
только одновременно.    
\end{propos}
Доказательство предложения \ref{prb2}, легко следующее   из определений $\ell^{\rh}$ и $\ell^{\lh}$ в \eqref{df:dDlm+}--\eqref{df:dDlm-} и условий Линделёфа \eqref{con:LpZR}--\eqref{con:LpZ}, опускаем. 
 
{\it Выметание\/} рода $01$ распределения зарядов $\nu$  {\it на  замкнутую вертикальную полосу\/}  $\overline \strip_b$ ширины $2b\geq 0$ из \eqref{{strip}c} опишем  в пять шагов [b\ref{i1}]--[b\ref{i5}], применяя каждый шаг к распределению зарядов, полученному на предыдущем шаге:
\begin{enumerate}[{[b1]}]
\item\label{i1} $(-b)$-сдвиг $\nu_{\vec{-b}}$ распределения зарядов $\nu$;
\item\label{i2} выметание $\nu_{\vec{-b}}^{\bal^{01}_{\rh}}$  рода  $01$ из правой полуплоскости  $\CC_{ \rh}$ на $ \CC_{\overline \lh}$;
\item\label{i3} $2b$-сдвиг $\Bigl(\nu_{\vec{-b}}^{\bal^{01}_{\rh}}\Bigr)_{\vec{2b}}$
 распределения зарядов $\nu_{\vec{-b}}^{\bal^{01}_{\rh}}$;
\item\label{i4} выметание $\Bigl(\nu_{\vec{-b}}^{\bal^{01}_{\rh}}\Bigr)_{\vec{2b}}^{{\bal^{01}_{\lh}}}$
рода $01$ из левой полуплоскости $ \CC_{\lh}$ на $\CC_{\overline \rh}$;
\item\label{i5} $(-b)$-сдвиг $\biggl(\Bigl(\nu_{\vec{-b}}^{\bal^{01}_{\rh}}\Bigr)_{\vec{2b}}^{{\bal^{01}_{\lh}}}\biggr)_{\vec{-b}}$ распределения зарядов $\Bigl(\nu_{\vec{-b}}^{\bal^{01}_{\rh}}\Bigr)_{\vec{2b}}^{{\bal^{01}_{\lh}}}$. 
\end{enumerate} 
Полученное  в [b\ref{i5}]  распределение зарядов для краткости обозначаем как  
\begin{equation}\label{Balb1}
\nu^{\Bal^{01}_b}:=\biggl(\Bigl(\nu_{\vec{-b}}^{\bal^{01}_{\rh}}\Bigr)_{\vec{2b}}^{{\bal^{01}_{\lh}}}\biggr)_{\vec{-b}}
\end{equation}
и называем {\it выметанием рода\/ $01$ на\/ $\overline \strip_b$ распределения зарядов\/} $\nu$, если  шаги  
[b\ref{i2}] и [b\ref{i4}] реализуемы.  Для последнего  по теореме  \ref{theoB1} достаточно, чтобы распределение зарядов $\nu$ было 
из класса сходимости при порядке $p\overset{\eqref{sufc}}{=}2$. 

\begin{remark}\label{remB01}
По замечанию \ref{rem0nu} при $\pm b\notin \supp \nu$
или, более общ\'о,  при 
\begin{equation}\label{Reznub}
\int_{b+r_0\DD} \Re^+\frac{1}{z-b} \dd |\nu|(z)+ 
\int_{-b+r_0\DD} \Re^-\frac{1}{z+b} \dd |\nu|(z)
<+\infty  
\end{equation}
для некоторого $r_0>0$ можно  в [b\ref{i2}] и [b\ref{i4}], а в итоге и в  \eqref{Balb1}  обойтись выметанием рода $1$, результат чего  в  \eqref{Balb1} будем обозначать   через $\nu^{\Bal^1_b}$. В частности, \eqref{Reznub} выполнено, если 
для некоторого числа $r_0>0$ имеем равенство
\begin{equation}\label{nur0+b}
|\nu|\bigl(b+r_0\DD\cap \CC_{ \rh}\bigr)+|\nu|\bigl(-b+r_0\DD\cap \CC_{ \lh}\bigr)=0.
\end{equation}

\end{remark}

\subsection{Сдвиги и выметание $\delta$-субгармонической функции}
Для точек $w\in \CC$ аналогично $w$-сдвигу \eqref{vecz} распределений зарядов определяем {\it $w$-сдвиг $u_{\vec{w}}$ 
функции\/} $u$ на $\CC$, задаваемый как 
\begin{equation}\label{veczu}
u_{\vec{w}}\colon z\underset{z\in \CC}{\longmapsto} u(z-w) .
\end{equation}
При $w$-сдвиге  \eqref{veczu}  $\delta$-субгармонической на $\CC$ функции $\mathcal U\not\equiv \pm \infty$ она остаётся такой же,  
а распределение её зарядов Рисса претерпевает $w$-сдвиг
\begin{equation}\label{zsh}
\frac{1}{2\pi}{\bigtriangleup}(\mathcal U_{\vec{w}})\overset{\eqref{vecz}}{=}
\Bigl(\frac{1}{2\pi}{\bigtriangleup}\mathcal U\Bigr)_{\vec{w}}.
\end{equation}  
Для $\delta$-субгармонической функции  $\mathcal U\not\equiv \pm\infty$ и $b\in \RR^+$ каждую  {\it $\delta$-субгар\-м\-о\-н\-и\-ч\-е\-с\-к\-ую  функцию\/} $\mathcal U^{\Bal_b}$,  равную  функции $\mathcal U$ на  вертикальной полосе $\overline \strip_b$ ширины $2b$ из \eqref{{strip}c} вне некоторого полярного множества   и
в то же время    гармоническую  на    $\CC\setminus \overline \strip_b$, называем 
{\it  $\delta$-субгармоническим выметанием на $\overline \strip_b$ функции\/  $\mathcal U$}. 

\begin{theorem}\label{Balvs}
 Пусть $b\in \RR^+$ и  $\delta$-суб\-г\-а\-р\-м\-о\-н\-и\-ч\-е\-с\-к\-ая  функция $\mathcal U\not\equiv \pm \infty$  c распределением зарядов
 Рисса  \eqref{Delta} ко\-н\-е\-ч\-н\-ой верхней плотности 
представима  в виде разности субгармонических функций не более чем первого порядка.
 Тогда существует $\delta$-субгармоническое  выметание\/ $\mathcal U^{\Bal_b}$ на $\overline \strip_b$
функции $\mathcal U$ c распределением зарядов Рисса 
\begin{equation}\label{112s}
\frac{1}{2\pi}{\bigtriangleup}(\mathcal U^{\Bal_b})\overset{\eqref{Balb1}}{=}\varDelta_\mathcal U^{\Bal^{01}_b},
\end{equation}  
представимое вне некоторого полярного множества  в виде разности 
\eqref{reprvB} субгармонических функций $u_{\pm}\not\equiv -\infty$ не более чем первого порядка. 
\end{theorem}
\begin{proof} Построение функции $\mathcal U^{\Bal_b}$ проводится в пять последовательных шагов
[B\ref{i1}]--[B\ref{i5}], применяя каждый шаг к $\delta$-субгармонической функции, полученной на предыдущем шаге:
\begin{enumerate}[{[B1]}]
\item\label{Bb1} $(-b)$-сдвиг $\mathcal U_{\vec{-b}}$ функции $\mathcal U$
с равенством $\frac{1}{2\pi}{\bigtriangleup}(\mathcal U_{\vec{-b}})\overset{\eqref{zsh}}{=}
(\varDelta_\mathcal U)_{\vec{-b}}$ и с очевидным сохранением для $\mathcal U_{\vec{-b}}$  представления в виде разности субгармонических функций порядка не больше  $1$;
\item\label{Bb2} $\delta$-субгармоническое выметание $\mathcal U_{\vec{-b}}^{\Bal_{\rh}}$ функции  
$\mathcal U_{\vec{-b}}$ из   $\CC_{ \rh}$ на $ \CC_{\overline \lh}$ в рамках теоремы 
\ref{Balv} с учётом равенства \eqref{112} в форме 
$\frac{1}{2\pi}{\bigtriangleup}\bigl({\mathcal U}^{\Bal_{\rh}}_{\vec{-b}}\bigr)
\overset{\eqref{112}}{=}(\varDelta_\mathcal U)_{\vec{-b}}^{\bal^{01}_{\rh}}$
и с представлением этого $\delta$-субгармонического выметания в виде разности субгармонических функций порядка не больше $1$;
\item\label{Bb3} $2b$-сдвиг $\bigl(\mathcal U_{\vec{-b}}^{\Bal_{\rh}}\bigr)_{\vec{2b}}$
функции $\mathcal U_{\vec{-b}}^{\Bal_{\rh}}$ и 
$
\frac{1}{2\pi}{\bigtriangleup}\bigl(\mathcal U_{\vec{-b}}^{\Bal_{\rh}}\bigr)_{\vec{2b}}
\overset{\eqref{zsh}}{=}\bigl((\varDelta_\mathcal U)_{\vec{-b}}^{\bal^{01}_{\rh}}\bigr)_{\vec{2b}}
$
с очевидным сохранением для $(\mathcal U_{\vec{-b}}^{\Bal_{\rh}})_{\vec{2b}}$ представления в виде разности субгармонических функций порядка не больше $1$;
\item\label{Bb4} зеркальная симметрия относительно $i\RR$ и применение выметания из $\CC_{ \rh}$ на $\CC_{\overline \lh}$ 
с обратной зеркальной симметрией позволяет определить $\delta$-суб\-г\-а\-р\-м\-о\-н\-и\-ч\-е\-с\-к\-ое  выметание  из $\CC_{ \lh}$ на $\CC_{\overline \rh}$ c  естественной верхней индексацией вида $^{\Bal_{\lh}}$,  а применение такого выметания 
в рамках теоремы  \ref{Balv} к функции $\bigl(\mathcal U_{\vec{-b}}^{\Bal_{\rh}}\bigr)_{\vec{2b}}$
с учётом равенства \eqref{112}  даёт 
$$
\frac{1}{2\pi}{\bigtriangleup}\bigl(\mathcal U_{\vec{-b}}^{\Bal_{\rh}}\bigr)_{\vec{2b}}^{\Bal_{\lh}}
\overset{\eqref{112}}{=}\bigl((\varDelta_\mathcal U)_{\vec{-b}}^{\bal^{01}_{\rh}}\bigr)_{\vec{2b}}^{\bal^{01}_{\lh}}
$$
вместе  с представлением для $\bigl(\mathcal U_{\vec{-b}}^{\Bal_{\rh}}\bigr)_{\vec{2b}}^{\Bal_{\lh}}$
в виде разности субгармонических функций порядка не больше $1$;
\item\label{Bb5} $(-b)$-сдвиг $\Bigl(\bigl(\mathcal U_{\vec{-b}}^{\Bal_{\rh}}\bigr)_{\vec{2b}}^{\Bal_{\lh}}\Bigr)_{\vec{-b}}$ функции
 $\bigl(\mathcal U_{\vec{-b}}^{\Bal_{\rh}}\bigr)_{\vec{2b}}^{\Bal_{\lh}}$ с равенством 
\begin{equation}\label{DU}
\frac{1}{2\pi}{\bigtriangleup}\Bigl(\bigl(\mathcal U_{\vec{-b}}^{\Bal_{\rh}}\bigr)_{\vec{2b}}^{\Bal_{\lh}}\Bigr)_{\vec{-b}}
\overset{\eqref{zsh}}{=}\biggl(\Bigl((\varDelta_\mathcal U)_{\vec{-b}}^{\bal^{01}_{\rh}}\Bigr)_{\vec{2b}}^{\bal^{01}_{\lh}}\biggr)_{\vec{-b}}
\overset{\eqref{Balb1}}{=}
\varDelta_{\mathcal U}^{\Bal^{01}_b}
\end{equation}
 где по построению $\Bigl(\bigl(\mathcal U_{\vec{-b}}^{\Bal_{\rh}}\bigr)_{\vec{2b}}^{\Bal_{\lh}}\Bigr)_{\vec{-b}}
={\mathcal U}^{\Bal_b}$ --- некоторое  $\delta$-субгармон\-и\-ч\-е\-с\-к\-ое выметание $\delta$-субгармонической функции $\mathcal U$ на $\overline \strip_b$, допускающее представление в виде разности субгармонических функций порядка не больше $1$, а \eqref{DU} и есть равенство \eqref{112s}. 
\end{enumerate}
\end{proof}

\begin{remark}\label{rem010}
Если в условиях теоремы   \ref{Balvs}
 функция $\mathcal U$  гармоническая в двух открытых полукругах $b+r_0\DD\cap \CC_{ \rh}$ и $-b+r_0\DD\cap \CC_{ \lh}$
для некоторого $r_0>0$, то по варианту  \eqref{nur0+b} замечания \ref{remB01}   правую часть 
в \eqref{112s} можно заменить на выметание  $\varDelta_\mathcal U^{\Bal^{1}_b}$ рода $1$ на $\overline \strip_{b}$, т.е. 
 \begin{equation}\label{112s0}
\frac{1}{2\pi}{\bigtriangleup}\bigl(\mathcal U^{\Bal_b}\bigr)\overset{\eqref{112s}}{=}\varDelta_\mathcal U^{\Bal^{1}_b}.
\end{equation}  
\end{remark}

\section{Выметание субгармонической функции  конечного типа\\
 на объединение вертикальной полосы и вещественной оси}\label{bRu}

\begin{theorem}\label{prop4_1}
Пусть  $M\not\equiv -\infty$ --- субгармоническая функция конечного типа с распределением масс  Рисса $\varDelta_M$.
Тогда для любого  $s\in \RR^+$ существуют субгармоническая   функция $M_{\RR}$ конечного типа 
с распределением масс Рисса  $\varDelta_{M_{\RR}}=\frac{1}{2\pi}{\bigtriangleup}M_{\RR}$, сосредоточенным на $\RR\setminus [-s,s]$, 
со свойством 
\begin{equation}\label{MMR}
\sup_{1\leq r<R<+\infty}\Bigl|\ell_{\varDelta_M}(r,R) -\ell_{\varDelta_{M_{\RR}}}(r,R)\Bigr|<+\infty, 
\end{equation} 
и субгармоническая функция $M_{s}$  конечного типа с носителем  $\supp \varDelta_{M_{s}}\subset \overline \strip_s$ 
 распределения масс Рисса $\varDelta_{M_s}=\frac{1}{2\pi}{\bigtriangleup}M_s$, для которых  
\begin{align}
M(x) &\equiv M_{\RR}(x)+M_{s}(x)\quad\text{при всех  $x\in \RR$},
\label{{MRab}r}
\\
M(z) &\equiv M_{\RR}(z)+M_{s}(z)\quad\text{при всех  $z\in \overline \strip_s$}, 
\label{{MRab}i}\\
M(z) &\leq M_{\RR}(z)+M_{s}(z)\quad\text{при всех  $z\in \CC$}. 
\label{{MRab}leq}
\end{align}
\end{theorem}
\begin{proof}
Ключевую роль будет играть довольно общая 
\begin{lemma}[{\rm \cite[основная теорема]{Kha91},  
\cite[предложение 2.1]{Kha01l},  \cite[теорема 8]{KhaShm19}}]\label{BR} 
Если для  $p\in \RR^+$ и замкнутой системы лучей $S$ на $\CC$ с одной общей вершиной 
раствор любого открытого угла, дополнительного к $S$, т.е. связной компоненты в  $\CC\setminus  S$,  строго меньше, чем  $\pi /p$,
то для любой субгармонической функции $u\not\equiv -\infty$ конечного типа 
при порядке $p$ существует субгармоническая функция $u^{\bal}\geq u$ на $\CC$ конечного типа 
при порядке $p$,  равная функции $u$ на каждом луче из $S$ и гармоническая в каждом дополнительном к  $S$ угле. 
\end{lemma}
Функция $u^{\bal}$ выше  ---  {\it выметание функции $u$ из открытого множества\/ $\CC\setminus  S$ на систему лучей\/ $S$.} При этом система лучей $S$ одновременно рассматривается и как замкнутое точечное множество в $\CC$ всех точек, лежащих на лучах из $S$.

Здесь применяем лемму \ref{BR}  при $p:=1$  в два  этапа. На первом этапе ---
 к функции $M$ в роли $u$  и бесконечной замкнутой системе  лучей $S_s^-$, 
состоящей из объединения луча $[s,+\infty)\subset \RR$ со всеми лучами из  замкнутой полуплоскости 
$\bigl\{z\in \CC\bigm| \Re z\leq s\bigr\}$ с общей вершиной в точке $s\in \RR$ и с  двумя дополнительными открытыми прямыми углами  раствора $\pi/2$.
 К полученной на первом этапе субгармонической функции  $M^{\bal}\geq M$ конечного типа на $\CC$, гармонической   на $\CC\setminus S_s^-$ и равной $M$ на $S_s^-$, применим  лемму \ref{BR} с функцией  $M^{\bal}$ 
в роли $u$  и бесконечной замкнутой системе  лучей $S_{-s}^+$, 
состоящей из объединения луча $(-\infty, -s]\subset \RR$ со всеми лучами из  замкнутой полуплоскости 
$\bigl\{z\in \CC\bigm| \Re z\geq -s\bigr\}$ с общей вершиной в точке $-s\in \RR$ и тоже с двумя дополнительными открытыми прямыми углами. 
В результате получим субгармоническую функцию 
$(M^{\bal})^{\bal}\geq M^{\bal} \geq M$ конечного типа на $\CC$, которая совпадает с $M$ на пересечении 
$S_s^-\cap S_{-s}^+=\RR\cup \overline\strip_s$ и гармоническая вне этого пересечения.  
По первой части  теоремы Вей\-е\-р\-ш\-т\-р\-а\-с\-са\,--\,Ада\-м\-а\-ра\,--\,Ли\-н\-д\-е\-л\-ё\-фа\,--\,Бре\-ло 
распределение масс  Рисса   $\frac{1}{2\pi}{\bigtriangleup}(M^{\bal})^{\bal}$ конечной верхней плотности, удовлетворяющее условию Линделёфа  и  сосредоточенное  на  $\RR\cup \overline \strip_s$, может быть разбито на сумму  сужения $\mu_{s}:=
\frac{1}{2\pi}{\bigtriangleup}(M^{\bal})^{\bal}{\lfloor}_{\overline \strip_s}$ на замкнутую полосу $\overline \strip_s$ и оставшуюся часть $\mu_{\RR}$, сосредоточенную на $\RR\setminus [-s,s]$. Каждое  из этих двух  распределений масс  $\mu_{s}$ и $\mu_{\RR}$, очевидно, конечной верхней плотности, а  их сумма $\mu_s+\mu_{\RR}=\frac{1}{2\pi}{\bigtriangleup}(M^{\bal})^{\bal}$ удовлетворяет условию Линделёфа.  Поэтому ввиду  сосредоточенности  $\mu_{\RR}$ на $\RR$  распределение масс $\mu_s$ удовлетворяет  $i\RR$-условию Линделёфа 
\begin{equation}\label{iRmab}
\sup_{r\geq 1}\biggl|\int_{1< |z|\leq r} \Im \frac{1}{ z}\dd \mu_{s}(z)\biggr|\overset{\eqref{con:LpZiR}}{<}+\infty. 
\end{equation}
Распределение масс  $\mu_{s}$ удовлетворяет и $\RR$-условие Линделёфа \eqref{con:LpZR}, так как   
\begin{multline}\label{LRab}
\sup_{r>1}\biggl|\int_{1< |z|\leq r} \Re \frac{1}{ z}\dd \mu_{s}(z)\biggr|
\leq \int_{|z|> 1}  \biggl|\Re \frac{1}{ z}\biggr|\dd \mu_{s}(z)\\
\overset{\eqref{df:nup}}{\leq} |s|\int_1^{+\infty}\frac{\dd \mu^{\rad}_{s}(t)}{t^2}
\leq  2|s|\int_1^{+\infty}\frac{\mu^{\rad}_s(t)}{t^3}\dd t
<+\infty
\end{multline}
для распределения масс $\mu_{s}$ конечной верхней плотности.  
Следовательно,  распределение масс $\mu_{s}$  удовлетворяет условие Линделёфа \eqref{con:LpZ}, как  и  распределение масс $\mu_\RR$, являющееся  разностью двух распределений масс, удовлетворяющих  условию Линделёфа \eqref{con:LpZ}. По  
 второй части   теоремы Вей\-е\-р\-ш\-т\-р\-а\-с\-са\,--\,Ада\-м\-а\-ра\,--\,Ли\-н\-д\-е\-л\-ё\-фа\,--\,Бре\-ло 
существует субгармоническая функция $M_{\RR}$ конечного типа с мерой Рисса $\frac{1}{2\pi}{\bigtriangleup}M_{\RR}=\mu_{\RR}$, а субгармоническая функция $M_{s}:=(M^{\bal})^{\bal} -M_{\RR}$ с мерой Рисса 
$\frac{1}{2\pi}{\bigtriangleup}M_{s}=\mu_{s}$ также конечного типа. 
За исключением пока свойства \eqref{MMR}, по построениям 
все остальные требования предложения \ref{prop4_1} к функциям $M_\RR$ и $M_{s}$, включая \eqref{{MRab}r}--\eqref{{MRab}leq}, выполнены. 

\begin{lemma}[{\cite[предложение 4.1, (4.19)]{KhaShmAbd20}}]\label{lemJl} 
Для  любой субгармонической функции $u\not\equiv -\infty$ конечного типа 
с распределением масс  Рисса $\varDelta_u$ в обозначении 
\begin{equation}\label{JiR}
J_{i\RR}(r,R; u):=\frac{1}{2\pi}\int_r^R \frac{u(-iy)+u(iy)}{y^2} \dd y, \quad 0<r<R\leq +\infty
\end{equation}
при любом $r_0\in \RR^+\setminus 0$ выполнены  соотношения 
\begin{align}
\sup_{r_0\leq r<R<+\infty} &\Bigl|J_{i\RR}(r,R;u)-\ell_{\varDelta_u}^{\rh}(r,R)\Bigr|
<+\infty,
\label{{Jll}r}
 \\
\sup_{r_0\leq r<R<+\infty} &\Bigl|J_{i\RR}(r,R;u)-\ell_{\varDelta_u}^{\lh}(r,R)\Bigr|
<+\infty,
\label{{Jll}l}
\\
\sup_{r_0\leq r<R<+\infty} &\Bigl|J_{i\RR}(r,R;u)-\ell_{\varDelta_u}(r,R)\Bigr|
<+\infty.
\label{{Jll}m} 
\end{align}
\end{lemma}

Соотношение \eqref{{Jll}r} этой леммы \ref{lemJl}, применённое  к $(-s)$-сдвигу
$M_{\vec{-s}}$ функции $M$ на $-s$, определённому в \eqref{veczu}, 
даёт соотношение 
\begin{equation}\label{lk1}
\sup_{1\leq r<R<+\infty}\Bigl|\ell^{\rh}_{\varDelta_{M_{\vec{-s}}}}(r,R)- 
J_{i\RR}(r,R;M_{\vec{-s}})\Bigr|<+\infty ,
\end{equation}
где $\varDelta_{M_{\vec{-s}}}$ --- это $(-s)$-сдвиг распределения масс из \eqref{vecz}.
Применение  к такому же сдвигу $(M_{\RR}+M_{s})_{\vec{-s}}$ даёт   соотношение
\begin{equation}\label{lk2}
\sup_{1\leq r<R<+\infty}\Bigl|J_{i\RR}\bigl(r,R;(M_{\RR}+M_{s})_{\vec{-s}}\bigr)
-\ell^{\rh}_{\varDelta_{(M_{\RR}+M_{s})_{\vec{-s}}}}(r,R) \Bigr|<+\infty .
\end{equation}
По  \eqref{{MRab}i} функции $(M_{\RR}+M_{s})_{\vec{-s}}$ и $M_{\vec{-s}}$
совпадают на $i\RR$ и, складывая   \eqref{lk1} с \eqref{lk2}, после оценок для суммы супремумов модулей получаем 
\begin{equation}\label{lk3}
\sup_{1\leq r<R<+\infty}\Bigl|\ell^{\rh}_{\varDelta_{M_{\vec{-s}}}}(r,R)- 
\ell^{\rh}_{\varDelta_{(M_{\RR}+M_{s})_{\vec{-s}}}} (r,R)
\Bigr|<+\infty. 
\end{equation}
Но   распределение масс Рисса $\varDelta_{(M_s)_{\vec{-s}}}$  сдвига $(M_s)_{\vec{-s}}$ функции $M_s$ сосредоточена на  $ \CC_{\overline \lh}$ и по определению правой логарифмической меры  
$$
\ell^{\rh}_{\varDelta_{(M_{\RR}+M_s)_{\vec{-s}}}}\overset {\eqref{df:dDlm+}}{=}\ell^{\rh}_{\varDelta_{(M_{\RR})_{\vec{-s}}}}.
$$
Отсюда согласно  \eqref{lk3} следует 
\begin{equation}\label{rop}
\sup_{1\leq r<R<+\infty}\Bigl|\ell^{\rh}_{\varDelta_{M_{\vec{-s}}}}(r,R)- 
\ell^{\rh}_{\varDelta_{(M_{\RR})_{\vec{-s}}}} (r,R)
\Bigr|<+\infty. 
\end{equation}
Ввиду конечности верхней плотности  распределений масс   $\varDelta_M$ 
из предложения \ref{prb2} следуют соотношения  
\begin{gather*}
\sup_{1\leq r<R<+\infty}\Bigl|\ell^{\rh}_{\varDelta_M}(r,R)- 
\ell^{\rh}_{\varDelta_{M_{\vec{-s}}}}(r,R)\Bigr|<+\infty, \\
\sup_{1\leq r<R<+\infty}\Bigl|\ell^{\rh}_{\varDelta_{(M_{\RR})_{\vec{-s}}}} (r,R)-
\ell^{\rh}_{\varDelta_{M_{\RR}}} (r,R)\Bigr|<+\infty. 
\end{gather*}
Cкладывая эти соотношения  с \eqref{rop} и  оценивая суммы супремумов модулей снизу через супермум модуля суммы,  получаем \eqref{MMR}.
\end{proof}

\end{fulltext}

\end{document}